\documentclass{birkjour}

 \newtheorem{theorem}{Theorem}
 \newtheorem*{theorem*}{Theorem}
 \newtheorem{corollary}{Corollary}
 
 \newtheorem{statement}{Statement}
 \newtheorem*{statement*}{Statement}
 \theoremstyle{definition}
 
 \theoremstyle{remark}
 \newtheorem{remark}{Remark}
 
 \numberwithin{equation}{section}

\begin{document}

	\title[Varignon's and Wittenbauer's parallelograms]
		{Varignon's and Wittenbauer's parallelograms}

	\author[Yuriy Zakharyan]{Yuriy Zakharyan}

	\address{%
		Leninskie Gory 1\\
		119234 Moscow\\
		Russian Federation\\
		ORCID: 0000-0003-0042-3372
	}

	\email{yuri.zakharyan@gmail.com}



	\date{June 1, 2018}

	\begin{abstract}
		In this paper the concept of homothetic parallelogram is introduced. This concept is a generalization of Varignon's and Wittenbauer's parallelograms of an arbitrary quadrangle, whose diagonals are not parallel. A formula for the area and perimeter of a homothetic parallelogram for the case when quadrangles are not crossed is obtained. The fact that homothetic parallelograms are similar to one another and are in perspective from diagonals intersection point is proved.
	\end{abstract}

	\maketitle
	\section{Introduction}
		We start with an arbitrary quadrangle $ABCD$, diagonals of which are not parallel. Let its diagonals intersect at point $O$. In this article we will focus on two related theorems.
		\begin{theorem*}[\textbf{Varignon}]
			\label{th:varignon}
			Midpoints of the sides of an arbitrary quadrangle form a parallelogram. It is called Varignon's parallelogram (Fig.1-a). ~[1, p.53]
		\end{theorem*}
		\begin{theorem*}[\textbf{Wittenbauer}]
			\label{th:wittenbauer}
			Let the sides of an arbitrary quadrangle be divided into three equal parts. Lines that join dividing points near its vertices form a parallelogram.
			It is called Wittenbauer's parallelogram (Fig.1-b). ~[2, p.216]
		\end{theorem*}
		\begin{center}
			\includegraphics[width=1\textwidth]{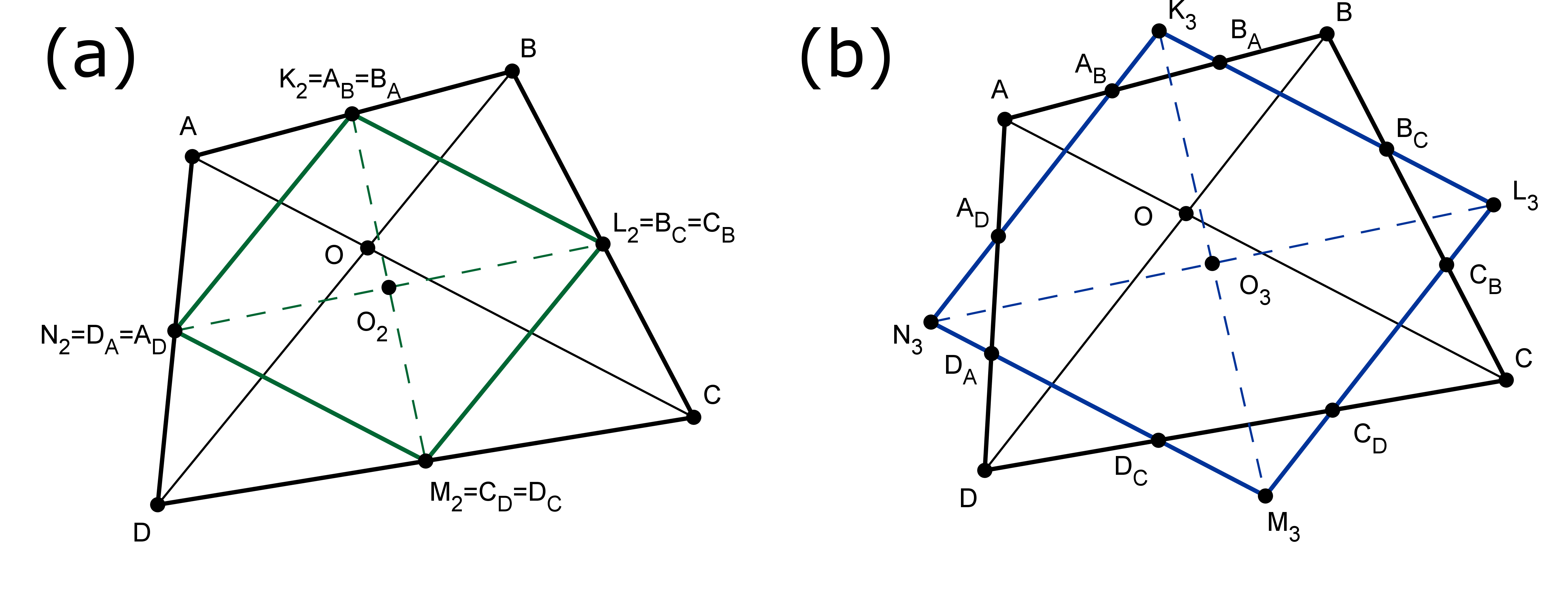}
			\newline
			Fig.1~Varignon's and Wittenbauer's parallelograms
		\end{center}
		It should be noticed that the quadrangle can be convex, re-entrant (Fig.2-a) or crossed (Fig.2-b). ~[1, p.52] 
		\begin{center}
			\includegraphics[width=1\textwidth]{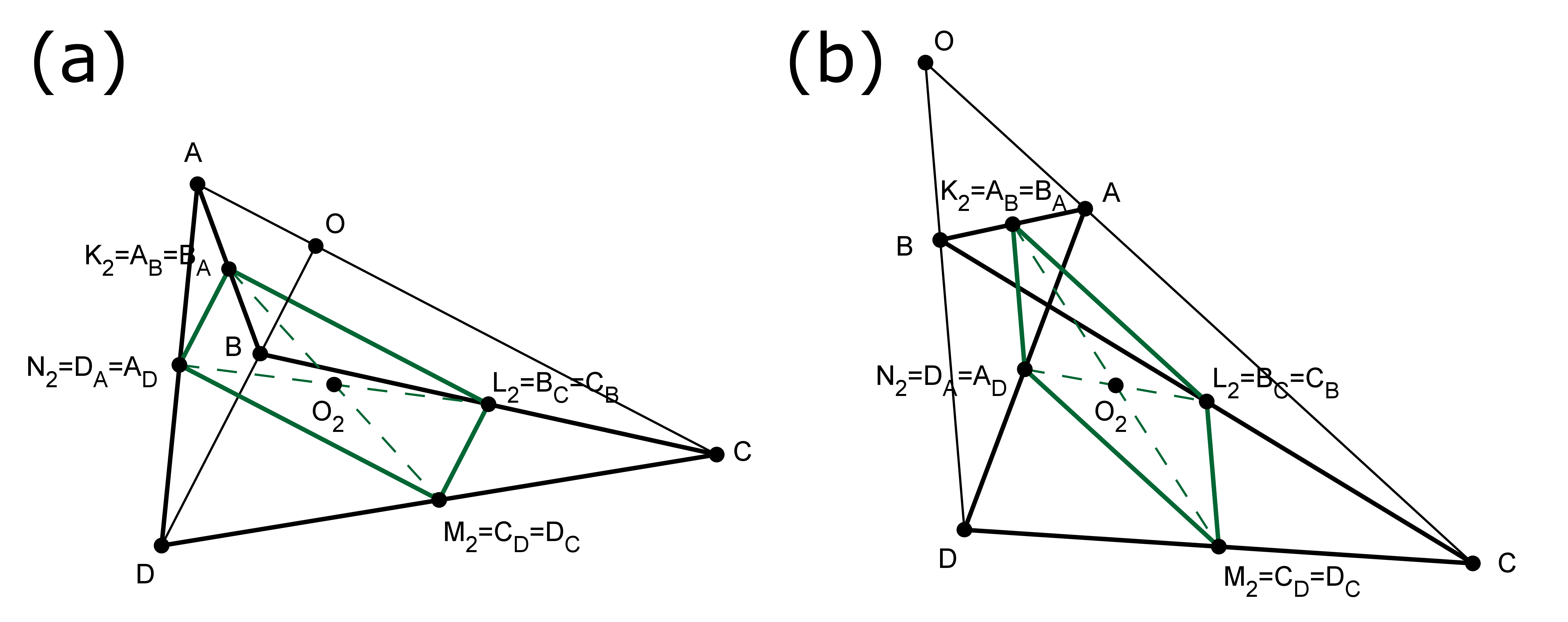}
			\newline
			Fig.2~Varignon's parallelogram of different quadrangles
		\end{center}
		There are a few statements related to Varignon's and Wittenbauer's parallelograms. 
		\begin{statement*}
			Let a quadrangle be convex or re-entrant and its area be $S_{ABCD}$. 
			The area of Varignon's parallelogram is $\frac{1}{2}S_{ABCD}$.
			The area of Wittenbauer's parallelogram is $\frac{8}{9}S_{ABCD}$. 
		\end{statement*}
		\noindent In fact, for a crossed quadrangle these areas are propotional to the difference between triangle areas with the same coefficients. 
		\begin{statement*}
			Varignon's and Wittenbauer's parallelograms are rectangles if and only if the quadrangle diagonals are perpendicular (Fig.3-a).  
			Varignon's and Wittenbauer's parallelograms are rhombuses if and only if the quadrangle diagonals are equal (Fig.3-b). 
		\end{statement*}
		\begin{center}
			\includegraphics[width=0.8\textwidth]{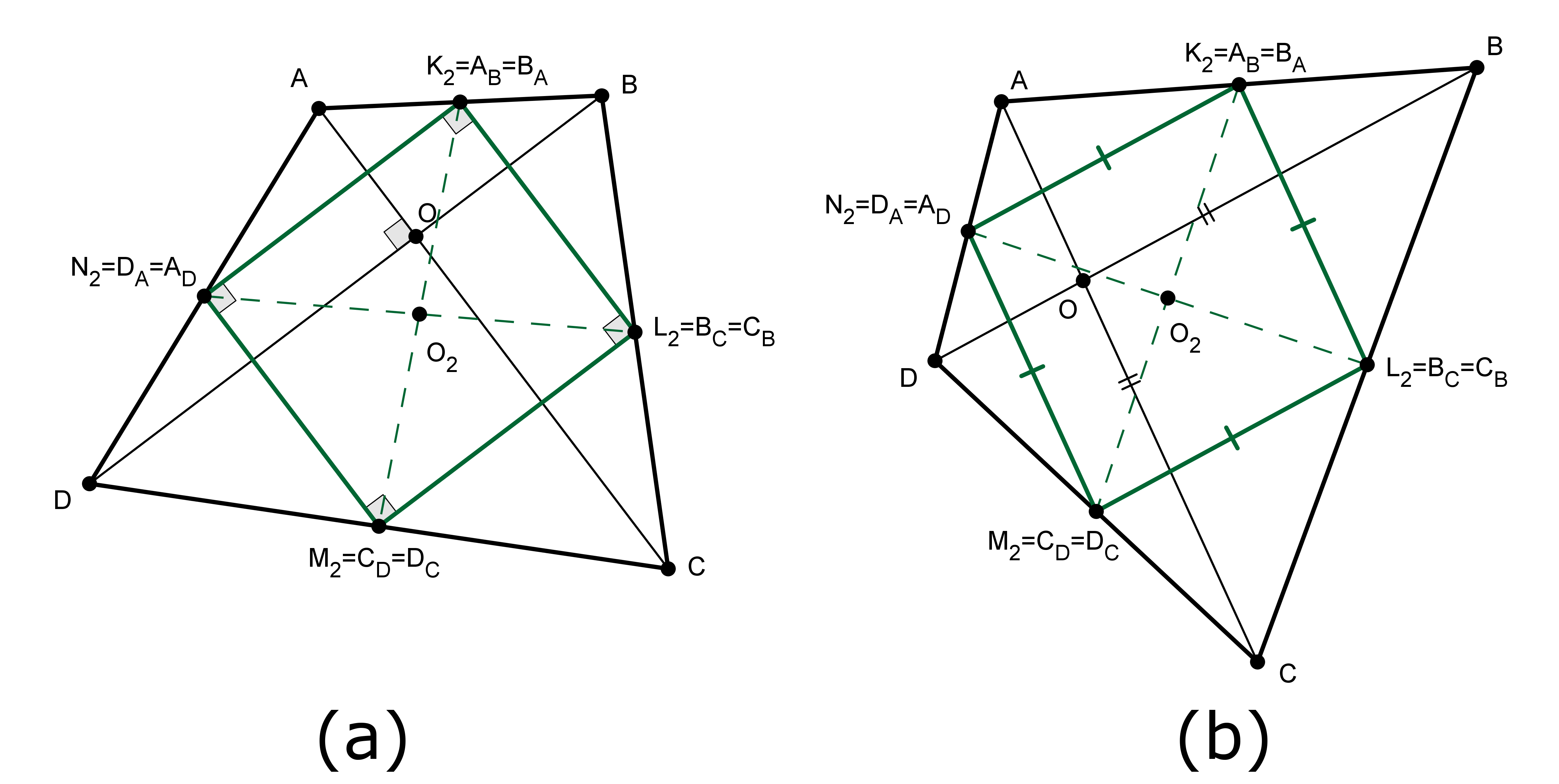}
			\newline
			Fig.3~Varignon's rectangle and rhombus
		\end{center}
	\section{Generalization}
		Varignon's and Wittenbauer's theorems are obviously related. In both theorems the quadrangle sides are divided into equal parts, dividing points are joined by lines and these lines form parallelograms.
		This similarity allows us to generalize these theorems. 
		\newline
		Let us divide the quadrangle sides into $n \in \mathbb{N}, n \geq 2$ equal parts. We can notice that the lines which join dividing points near the vertices form a parallelogram (Fig.4). 
		\begin{center}
			\includegraphics[width=0.4\textwidth]{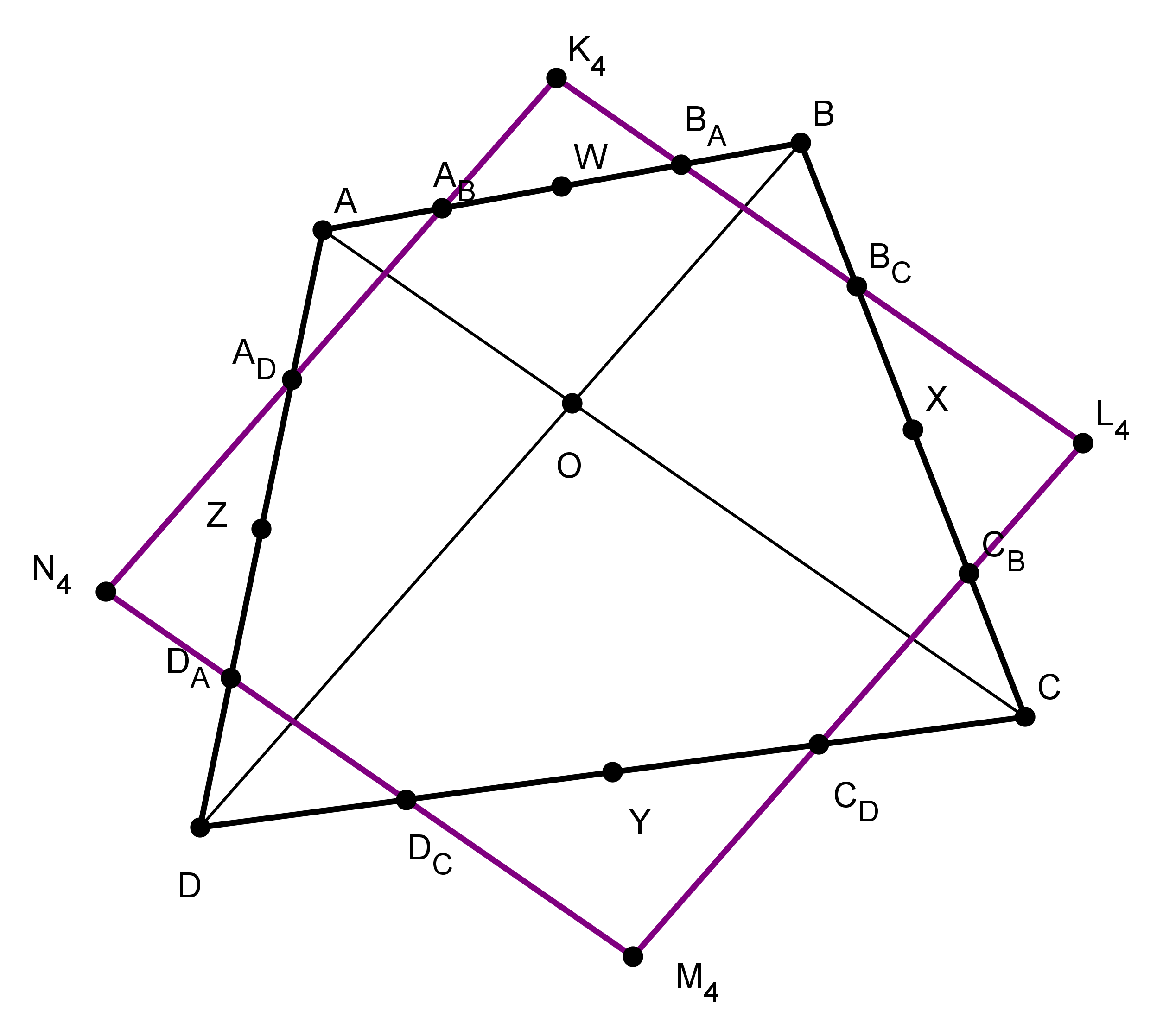}
			\newline
			Fig.4~Parallelogram, n=4
		\end{center}
		Let us notice that we join only dividing points near the vertices. And the distance between these points and respective vertices is proportional to the quadrangle sides.
		From this point of view these dividing points are homothetically transformed. ~[3, p.145]
		\newline
		Let $A_{B}^{\lambda}$ (or $A_{B}(\lambda)$) mean homothetic tranformation of the point $B$ with center $A$ and ration $\lambda$. Let us formulate the first theorem.
		\begin{theorem}
			\label{th:th1}
			$\forall \lambda \in \mathbb{R}$  lines 
			$A_{D}^{\lambda}A_{B}^{\lambda}$, $B_{A}^{\lambda}B_{C}^{\lambda}$, $C_{B}^{\lambda}C_{D}^{\lambda}$, $D_{C}^{\lambda}D_{A}^{\lambda}$ form 
			a parallelogram. Let us call it a homothetic parallelogram of the $ABCD$ with ratio $\lambda$.
		\end{theorem}
		\begin{center}
			\includegraphics[width=1\textwidth]{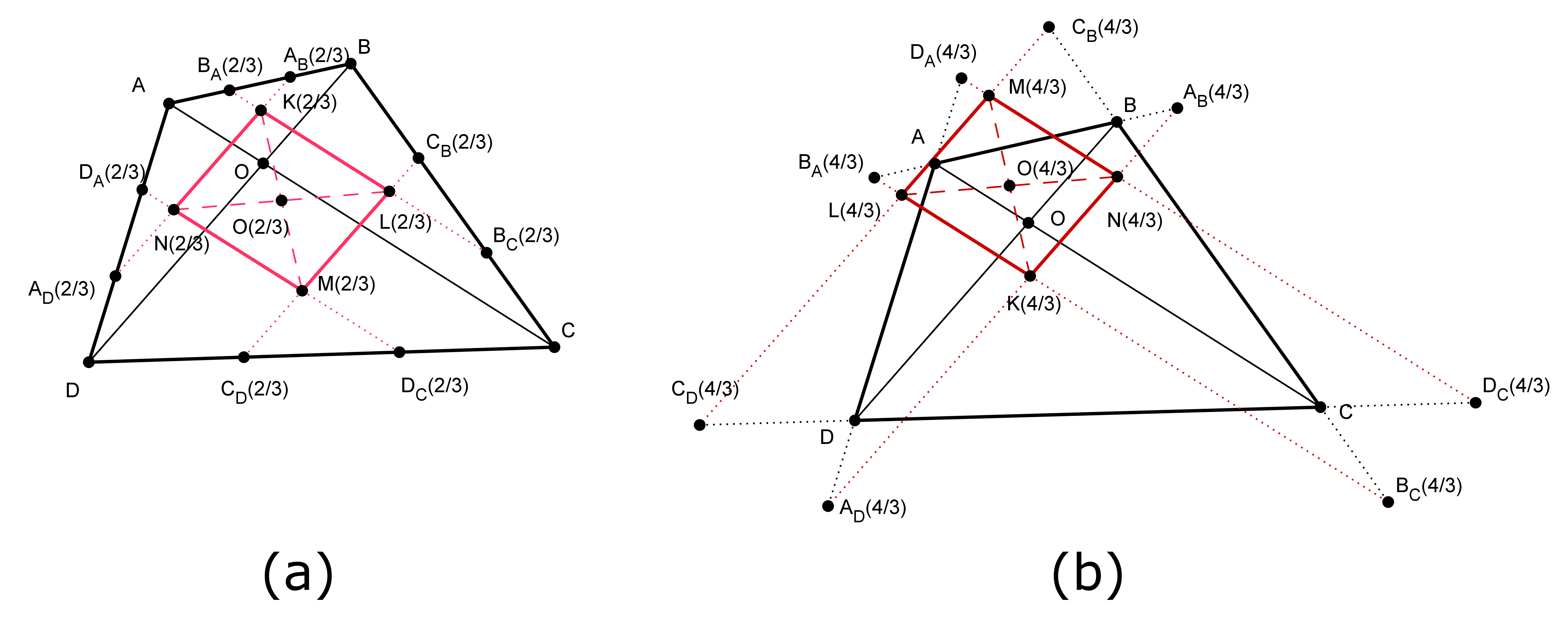}
			\newline
			Fig.5~Homothetic parallelograms, $\lambda=\frac{2}{3},\frac{4}{3}$
		\end{center}
		\begin{proof}
			From homothety it follows:
			\begin{equation*}
				\begin{aligned}
					\frac{\left |BB_{A}^{\lambda}\right |}{\left |BA\right |}&=\left|\lambda\right|\\
					\frac{\left |BB_{C}^{\lambda}\right |}{\left |BC\right |}&=\left|\lambda\right|
				\end{aligned}
			\end{equation*}
			Moreover, ${\Delta}ABC, {\Delta}{B_{A}^{\lambda}}B{B_{C}^{\lambda}}$ have a common angle. Therefore ${\Delta}ABC\sim{\Delta}{B_{A}^{\lambda}}B{B_{C}^{\lambda}}$ and $AC \mid \mid {B_{A}^{\lambda}}{B_{C}^{\lambda}}$. ~[2, p.8]
			\newline
			By analogy, 
			\begin{equation*}
				\begin{aligned}
					AC &\mid \mid {D_{A}^{\lambda}}{D_{C}^{\lambda}}\\
					BD &\mid \mid {A_{B}^{\lambda}}{A_{D}^{\lambda}}\\
					BD &\mid \mid {C_{B}^{\lambda}}{C_{D}^{\lambda}}
				\end{aligned}
			\end{equation*}
			Finally,
			\begin{equation*}
				\begin{aligned}
					{B_{A}^{\lambda}}{B_{C}^{\lambda}} &\mid \mid {D_{A}^{\lambda}}{D_{C}^{\lambda}}\\
					{A_{B}^{\lambda}}{A_{D}^{\lambda}} &\mid \mid {C_{B}^{\lambda}}{C_{D}^{\lambda}}
				\end{aligned}
			\end{equation*}
			\qedhere
		\end{proof}
		From the proof of Theorem~\ref{th:th1} we have
		\begin{corollary}
			\label{cor:cor1}
			The sides of a homothetic parallelogram are parallel to the quadrangle diagonals. 
			Also, homothetic parallelograms are similar to each other.
		\end{corollary}
		\begin{corollary}
			\label{cor:cor2}
			The homothetic parallelogram is a rectangle if and only if the quadrangle diagonals are perpendicular.  
			The homothetic parallelogram is a rhombus if and only if the quadrangle diagonals are equal. 
		\end{corollary}
		\begin{remark}
			Varignon's parallelogram is a homothetic parallelogram with ratio ${\lambda}=\frac{1}{2}$.
			Wittenbauer's parallelogram is a homothetic parallelogram with ratio ${\lambda}=\frac{1}{3}$.
			\newline
			The homothetic parallelogram with ratio $\lambda=1$ is the diagonals intersection point (Fig.6-a), the homothetic parallelogram with ratio $\lambda=0$ is a limiting parallelogram (Fig.6-b). 
		\end{remark}
		\begin{center}
			\includegraphics[width=0.8\textwidth]{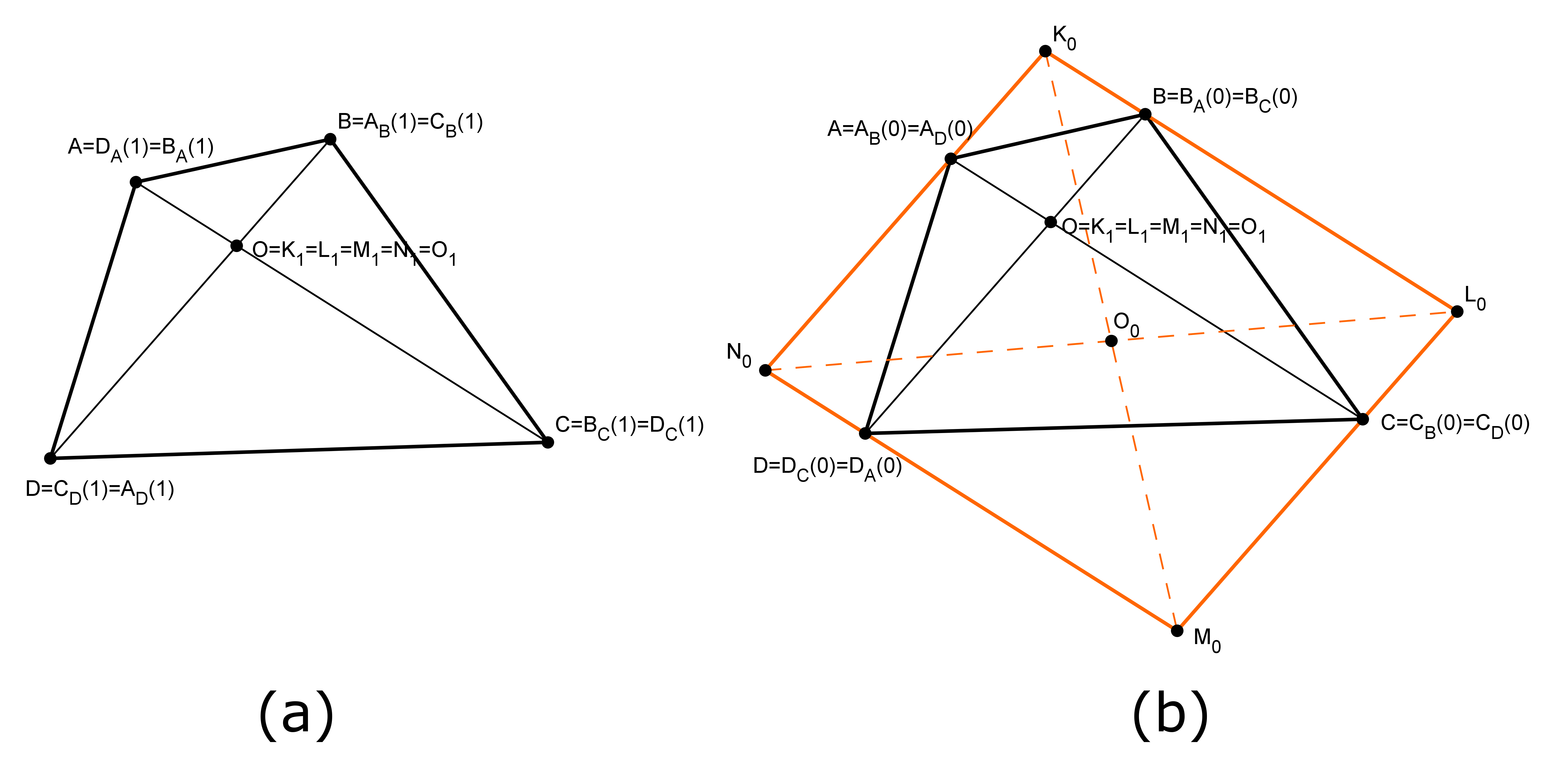}
			\newline
			Fig.6~Homothetic parallelograms, $\lambda=1,0$
		\end{center}
		\begin{remark}
			From this moment let a homothetic parallelogram with ratio $\lambda$ be $K^{\lambda}L^{\lambda}M^{\lambda}N^{\lambda}$, where
			\begin{equation*}
				\begin{aligned}
					K^{\lambda}&={A_{B}^{\lambda}}{A_{D}^{\lambda}}\cap{B_{A}^{\lambda}}{B_{C}^{\lambda}}\\
					L^{\lambda}&={B_{A}^{\lambda}}{B_{C}^{\lambda}}\cap{C_{B}^{\lambda}}{C_{D}^{\lambda}}\\
					M^{\lambda}&={C_{B}^{\lambda}}{C_{D}^{\lambda}}\cap{D_{A}^{\lambda}}{D_{C}^{\lambda}}\\
					N^{\lambda}&={D_{A}^{\lambda}}{D_{C}^{\lambda}}\cap{A_{B}^{\lambda}}{A_{D}^{\lambda}}	
				\end{aligned}
			\end{equation*}
		\end{remark}
		For homothetic parallelograms there is also a formula for the area. 
		\begin{statement}
			\label{st:st1}
			Let a quadrangle be convex or re-entrant and its area be $S_{ABCD}$. 
			The area of a homothetic parallelogram is $2{(\lambda-1)}^{2}S_{ABCD}$.
		\end{statement}
		\begin{proof}
			The proof depends on ratio $\lambda$. Also, it depends on the type of a quadrangle. It is based on the similatiry of the triangles and summation-subtraction of the areas. Let us prove the statement for a convex quadrangle, with ratio $\lambda<0$ (Fig.7). The proof for other cases is analogous.
			\begin{center}
				\includegraphics[width=0.6\textwidth]{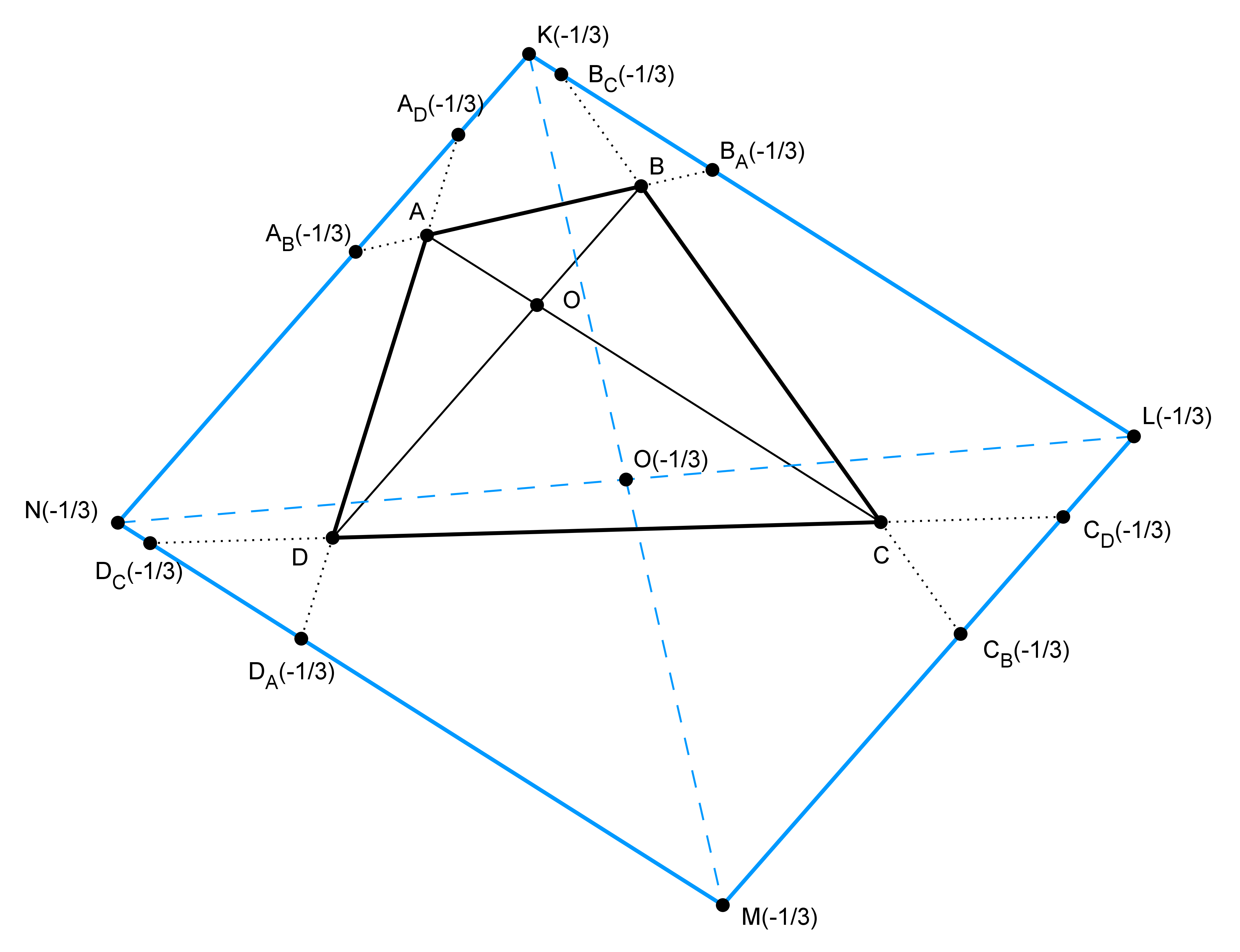}
				\newline
				Fig.7~Homothetic parallelogram, $\lambda=-\frac{1}{3}$
			\end{center}
			Considering $\lambda<0$, we obtain:
			\begin{equation*}
				\begin{aligned}
					S_{K^{\lambda}L^{\lambda}M^{\lambda}N^{\lambda}}&=S_{ABCD}+S_{K^{\lambda}{A_{B}^{\lambda}}{B_{A}^{\lambda}}}+S_{L^{\lambda}{B_{C}^{\lambda}}{C_{B}^{\lambda}}}+
					S_{M^{\lambda}{C_{D}^{\lambda}}{D_{C}^{\lambda}}}+S_{N^{\lambda}{D_{A}^{\lambda}}{A_{D}^{\lambda}}}-\\
					&-S_{A{A_{D}^{\lambda}}{A_{B}^{\lambda}}}-
					S_{B{B_{A}^{\lambda}}{B_{C}^{\lambda}}}-S_{C{C_{B}^{\lambda}}{C_{D}^{\lambda}}}-S_{D{D_{C}^{\lambda}}{D_{A}^{\lambda}}}=\\
					&=S_{ABCD}+(1-2\lambda)^2S_{OBA}+(1-2\lambda)^2S_{OCB}+\\
					&+(1-2\lambda)^2S_{ODC}+(1-2\lambda)^2S_{OAD}-\\
					&-{\lambda}^2S_{ADB}-{\lambda}^2S_{BAC}-{\lambda}^2S_{CBD}-{\lambda}^2S_{DCA}=\\
					&=\left ( {1+(1-2\lambda)^2-2{\lambda }^2} \right ) S_{ABCD}=\left ( {2-4\lambda+2{\lambda}^2} \right ) S_{ABCD}=\\ 
					&= 2{(\lambda-1)}^{2}S_{ABCD}
				\end{aligned}	
			\end{equation*}
			\qedhere
		\end{proof}
		\begin{remark}
			\label{rm:rm1}
			For the crossed quadrangle $ABCD$ there is an analogous formula. The area of homothetic parallelogram in this case is $2{(\lambda-1)}^{2}S$, where $S$ does not depend on $\lambda$. $S$ is the difference between the areas of the triangles. 
		\end{remark}
		\noindent Taking Corollary~\ref{cor:cor1}, Statement~\ref{st:st1} and Remark ~\ref{rm:rm1}, we have: 
		\begin{corollary}
			\label{cor:cor3}
			$\forall {\lambda}_{1}, {\lambda}_{2} \in \mathbb{R}:~\frac{p^{\lambda_1}}{p^{\lambda_2}}=\frac{\left | {\lambda_1 -1} \right |}{\left | {\lambda_2 -1} \right |}$, where $p^{\lambda_i}$ means the perimeter of a homothetic parallelogram with ratio $\lambda_i$. Here $\lambda_2$ can be formally equal to $1$.  
		\end{corollary}
	\section{Perspective}
		This section is related to the theory of perspective~[1, p.70]. 
		\begin{theorem}
			\label{th:th2}
			Homothetic parallelograms are in perspective from the diagonals intersection point $O$. (Fig.8) Moreover,
			\begin{equation*}
				\forall {\lambda}_{1}, {\lambda}_{2} \in \mathbb{R}:~\frac{\left |OK^{\lambda_1}\right |}{\left |OK^{\lambda_2}\right |}=\frac{\left |OL^{\lambda_1}\right |}{\left |OL^{\lambda_2}\right |}=
				\frac{\left |OM^{\lambda_1}\right |}{\left |OM^{\lambda_2}\right |}=\frac{\left |ON^{\lambda_1}\right |}{\left |ON^{\lambda_2}\right |}=\frac{\left | {\lambda_1 -1} \right |}{\left | {\lambda_2 -1} \right |}.
			\end{equation*}
			Here the denominator can be formally equal to $0$. 
		\end{theorem}
		\begin{center}
			\includegraphics[width=1\textwidth]{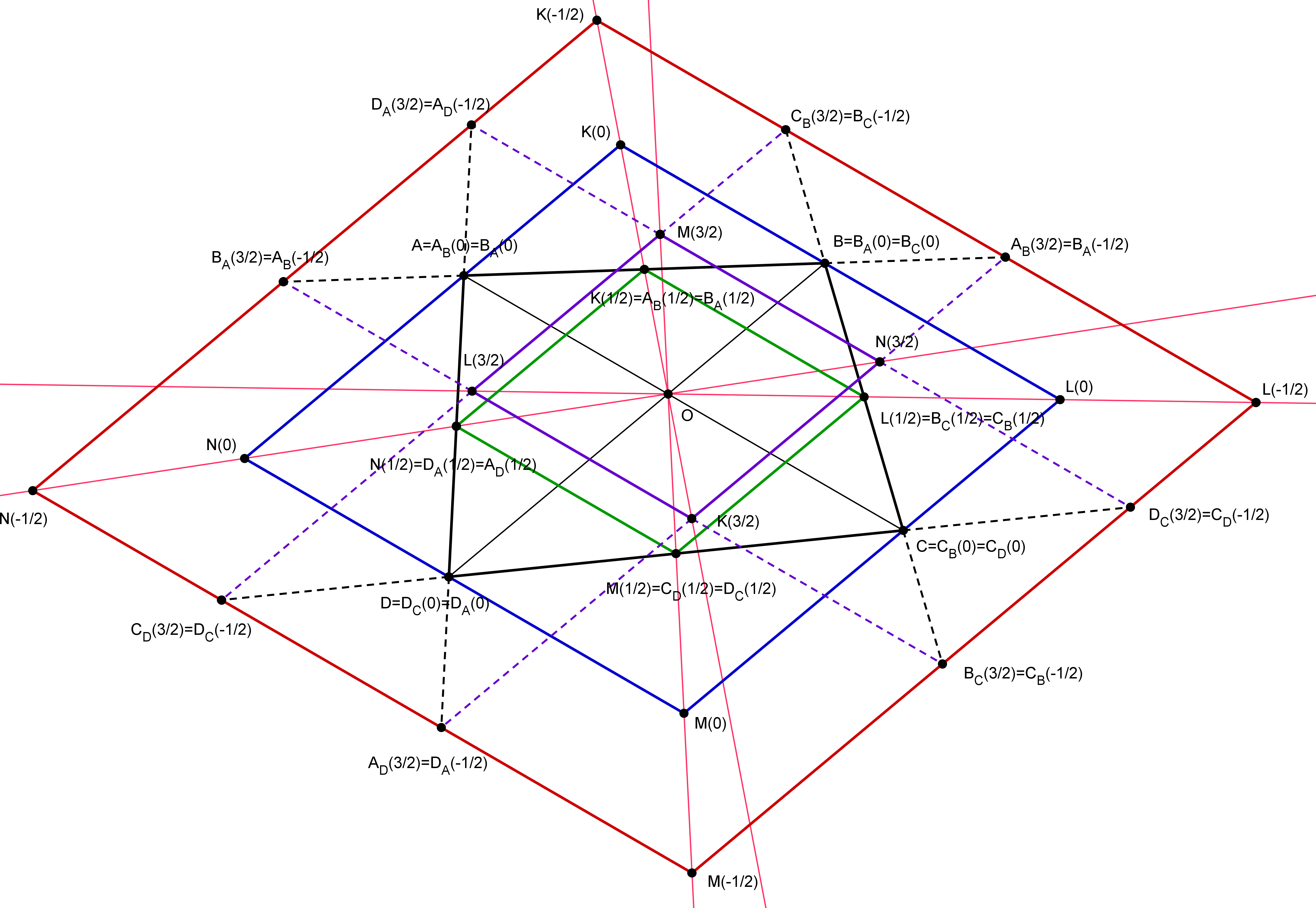}
			\newline
			Fig.8~Perspective parallelograms
		\end{center}
		\begin{proof}
			Let us prove that $\forall \lambda \in \mathbb{R}: K^{\lambda} \in OK^{0}$.
			\newline
			It is obvious that $K^{0}BOA$ is a parallelogram. Thus, $K^{\frac{1}{2}} \in OK^{0}$ as it is the midpoint of $AB$. 
			Also, it is obvious that $\forall \lambda \in \mathbb{R}:A_{B}^{\lambda}=B_{A}^{1-\lambda}$. 
			\newline
			Thus, the triangles ${\Delta}K^{\lambda}A_{B}^{\lambda}B_{A}^{\lambda}$ and ${\Delta}K^{1-\lambda}A_{B}^{1-\lambda}B_{A}^{1-\lambda}$ have a common side $A_{B}^{\lambda}B_{A}^{\lambda}=A_{B}^{1-\lambda}B_{A}^{1-\lambda}$ with the midpoint $K^{\frac{1}{2}}$. 
			\newline
			Moreover, other sides being parallel, we have  ${\Delta}K^{\lambda}A_{B}^{\lambda}B_{A}^{\lambda}={\Delta}K^{1-\lambda}A_{B}^{1-\lambda}B_{A}^{1-\lambda}$.
			As a result, $K^{\lambda}A_{B}^{\lambda}K^{1-\lambda}B_{A}^{\lambda}$ is a parallelogram and $K^{\lambda}, K^{\frac{1}{2}},K^{1-\lambda}$ are colinear. So, it is sufficient to prove it for the case 
			$\lambda<\frac{1}{2}$. Then $K^{0}, K^{\lambda}$ belong to one semiplane from $AB$.  It is obvious that ${\Delta}K^{\lambda}A_{B}^{\lambda}B_{A}^{\lambda}\sim{\Delta}K^0AB$ and their similar sides have a common midpoint. Thus, 
			\begin{equation*}
				{\angle}A_{B}^{\lambda}K^{\lambda}K^{\frac{1}{2}}={\angle}AK^{0}K^{\frac{1}{2}}
			\end{equation*}
			that means $K^{\lambda}, K^{0}, K^{\frac{1}{2}}$ are colinear. Finally, $K^{\lambda} \in OK^{0}$. For the other vertices the proof is analogous. 
			\newline
			With the proportionality theorem ~[3, p.116] we have
			\begin{equation*}
				\frac{\left |OK^{\lambda_1}\right |}{\left |OK^{\lambda_2}\right |}=\frac{\left |CB_{C}^{\lambda_1}\right |}{\left |CB_{C}^{\lambda_2}\right |}=\frac{\left |CC_{B}^{1-\lambda_1}\right |}{\left |CC_{B}^{1-\lambda_2}\right |}=\frac{\left | {1-\lambda_1} \right |}{\left | {1-\lambda_2} \right |}
			\end{equation*}
			We should notice that the denominator can be formally equal to zero. For the other vertices the proof is analogous.
			\qedhere
		\end{proof}
	
	\section{Conclusion}
		We have seen that Varignon's and Wittenbauer's parallelograms are related. Moreover, we can generalize them to homothetic parallelograms (Theorem \ref{th:th1}). The main properties of the parallelograms can also be generalized (Statement \ref{st:st1}, Corollaries \ref{cor:cor1}--\ref{cor:cor3}). It turned out that the homothetic parallelograms are in perspective from the diagonals intersection point (Theorem \ref{th:th2}) and the vertices ratio (the three colinear points ratio ~[4, p.29]) is predefined (Theorem \ref{th:th2}).
		\newline
		These results were presented at the European Union Contest for Young Scientists (EUCYS-2011) ~[5] and at the International Conference for Young Scientists (ICYS-2012) ~[6].
		\newline
		Later in the publication ~[7, pp.27-36] Romanian mathematician Kiss Sandor defined Wittenbauer-type parallelogram as a special case of the homothetic parallelogram (where $\lambda \in 
		\mathbb{Q}\cap\left(0, 1\right)$). For this special case he presented the proof of vertices colinearity, the formula for vertices coordinates, the parallelograms perimeter and area.

\end{document}